\documentstyle[12pt,epsf]{article}

\title{Five Definitions of Critical Point at Infinity}
\author{Alan H. Durfee}

\newcommand{\postscript}[2]
{\setlength{\epsfxsize}{#2\hsize}
\centerline{\epsfbox{#1}}}

\input amssym.def
\input amssym.tex
\def\real{{\Bbb R}}
\def\complex{{\Bbb C}}

\def\projective{{\Bbb P}}

\def\complexprojective{{\Bbb P}}

\def\lineinfinity{{\Bbb L}_\infty}
\def\clinf{{\Bbb L}_\infty}

\def\endofproof{$\square$}
\def\infinity{\infty}

\def\cinf{\complex \cup \{\infinity\}}

\newcounter{mycounter}[section]
\renewcommand{\themycounter}{\arabic{section}.\arabic{mycounter}}

{\medskip 
    \refstepcounter{mycounter}
    {\bf \noindent Theorem \themycounter. \ } \begin{em} }%
{\end{em} \medskip }

\newenvironment{proposition}%
{\medskip 
    \refstepcounter{mycounter}
    {\bf \noindent Proposition \themycounter. \ } \begin{em} }%
{\end{em} \medskip }

\newenvironment{lemma}%
{\medskip 
    \refstepcounter{mycounter}
    {\bf \noindent Lemma \themycounter. \ } \begin{em} }%
{\end{em} \medskip }

\newenvironment{corollary}%
{\medskip 
    \refstepcounter{mycounter}
    {\bf \noindent Corollary \themycounter. \ } \begin{em} }%
{\end{em} \medskip }

\newenvironment{remark}%
{\medskip 
    \refstepcounter{mycounter}
    {\bf \noindent Remark \themycounter. \ }}%
{\medskip }

\newenvironment{definition}%
{\medskip 
    \refstepcounter{mycounter}
    {\bf \noindent Definition \themycounter. \ }}%
{\medskip }

\newenvironment{example}%
{\medskip 
    \refstepcounter{mycounter}
    {\bf \noindent Example \themycounter. \ }}%
{\medskip }

{\medskip 
    \refstepcounter{mycounter}
    {\bf \noindent Problem \themycounter. \ }}%
{\medskip }

\newenvironment{xproof}%
{\medskip 
  \noindent
    {\bf Proof. \ }}%
{\endofproof \medskip }

\newcommand{\resetassertioncounter}{}

\begin{document}
%the next command causes all items in crptc.bib to be listed in the bibliography,
%referenced or not.
%\nocite{Bartolo-93,Bartolo-Cassou-Velasco-P93,Broughton-83,Broughton-88,Cassou-Ha-P1,Cassou-Ha-P2,Durfee-index,REU,Gaffney-P95,Ha-89-1,Ha-89-2,Ha-91,Ha-P91,Ha-P92,Ha-P90,Ha-Le,Krasinski,Libgober-93,Nemethi-Zaharia-90,Nemethi-Zaharia-92,Neumann-89,Neumann-Thanh-P92,Siersma-Tibar-P95,Suzuki-74}

\maketitle

\begin{abstract}
This survey paper discusses five equivalent ways of defining a ``critical point at
infinity'' for a complex polynomial of two variables. 
\end{abstract}

%\tableofcontents

%%%%%%%%%%%%%%%%%%%%%%%%%%%%%%%%%

\section{Introduction}

A proper smooth map without critical points from one manifold to
another is a locally trivial fibration by a well-known theorem of
Ehresmann.  On the other hand, a nonproper map without
critical points may not be a fibration.  This phenomenon occurs for
complex polynomials.  A simple example is provided by $f: \complex^2
\to \complex$ defined by the polynomial
$f(x,y) = y(xy-1)$.  This map has no critical
points, but the fiber over the origin is different from the other
fibers.  (In fact, the fiber over the origin is two rational curves,
one punctured at two points and the other at one point, whereas the
general fiber is a cubic curve, punctured at two points.)  One would like to
identify these ``critical values'' where the topology changes and their
corresponding ``critical points at infinity''.
We first review the history of this subject.

Let $f: \complex^{n} \to \complex$ be a complex polynomial.
There is a finite set $\Sigma \in \complex$ such that 
$$ f: \complex^n - f^{-1}(\Sigma) \to \complex - \Sigma $$
is a fibration.  
This is a form of Sard's theorem for polynomials; the set $\Sigma$ is finite
because it is algebraic.
For a proof, see \cite[Proposition 1]{Broughton-83}
(based on work of Verdier), \cite[Appendix A1]{Pham-83},
\cite[Theorem 1]{Ha-Le-84} or \cite{Ha-89-2}.
We let
$$\Sigma = \Sigma_{fin} \cup \Sigma_{\infinity}$$
where $\Sigma_{fin}$ is the set of critical values coming from critical points
in $\complex^n$, and $\Sigma_\infinity$ is the set of critical values
``coming from infinity''.
Of course these two sets may have nonempty intersection.

Broughton in \cite{Broughton-83,Broughton-88} calls the polynomial
$f$ {\em tame} if there is a $\delta > 0$ such that the set $\{ x: |
grad \, f(x) | \leq \delta \}$ is compact.  He proved that if $f$ is
tame, then $\Sigma_\infinity$ is empty.  

Thus if the gradient of a polynomial goes to zero along some path
going to infinity, then something bad may happen.
Topics surrounding the
gradient of the polynomial are treated in Section 4 of this paper.
The speed at which the gradient of $f$ goes to zero is measured by the
Lojasiewicz number at infinity; see \cite{Ha-90, Ha-P91, Ha-94,
Cassou-Ha-P92, Cassou-Ha-95}.

There followed many efforts in the case $n = 2$ to identify
the set $\Sigma_\infinity$ more precisely.
Suzuki \cite[Corollary 1]{Suzuki-74} provides an
estimate on the number of points in $\Sigma$.
In \cite{Ha-Le-84} it is shown that $c \in \Sigma$ if and
only if $\chi (f^{-1} (c) ) \neq \chi (f^{-1} (t) )$, where $f^{-1} (t)$ is a
generic fiber of $f$ and $\chi$ denotes Euler characteristic.
Further work on identifying $\Sigma_\infinity$ ca be found in
\cite{Ha-Nguyen-89, Ha-89-2, Nemethi-Zaharia-90, Nemethi-Zaharia-92,
Le-Oka-P93}.  

The homology and homotopy of the fibers of the polynomial $f$ were
also computed, leading to various numerical invariants which will be
discussed in the Section 2 of this paper.  Suzuki \cite[Proposition 2]{Suzuki-74} shows that
$$rank \, H_1(f^{-1}(t)) = \mu + \lambda$$
where $f^{-1}(t)$ is a generic fiber, $\mu$ is the sum of the Milnor numbers at the critical points of
$f$ in $\complex^2$, and $\lambda$ is the sum of all the ``jumps'' in
the Milnor numbers at infinity.  (In the terminology of Section 2, 
$\lambda = \sum \nu_{p,c}$, where the sum is over $c \in \complex$ and
$ p \in \lineinfinity$, and $\nu_{p,c}$ is the jump in the Milnor
number at the point $p \in \lineinfinity$ and value $c \in \complex$.) 
More on the topology of the fiber can be found in \cite{Bartolo-Cassou-Dimca-P96}.

The polynomial $f$ extends to a function on projective space
$\projective^2$ which is
well defined except at a finite number of points.  The points of
indeterminacy can be easily resolved, and the structure of the
resolution contains information about these points
\cite{Le-Weber-95, Le-Weber-P96}.
These topics are discussed in Section 3.

Other topics investigated (but not discussed in this paper) include
Newton diagrams \cite[Proposition 3.4]{Broughton-88},
\cite{Nemethi-Zaharia-90, Cassou-P96}, knots \cite{Neumann-89, Ha-91},
and the Jacobian conjecture \cite{Le-Weber-95}.

Papers in higher dimensions (that is, $n > 2$) include \cite{Parusinski-P95, Siersma-Tibar-95, Tibar-P96}.
Broughton \cite[Proposition 3.2]{Broughton-88} shows that the tame
polynomials form a dense constructible set in the set of polynomials
of a given degree; Cassou-Nogues \cite[Example V]{Cassou-P96} gives an example to
show that this set is not open in dimension $n=3$.

Although the scene in higher dimensions is not yet settled, the
situation in dimension two is now clear.  The purpose of this paper is
to collect together five definitions of ``critical point at infinity'' 
in this low-dimensional case
and prove that they are equivalent.
These definitions have appeared in the literature in some form
or other, usually in a global affine context;    
the purpose of this paper is to give these definitions and prove
their equivalence in a purely
local setting near a point on the line at infinity.
Many examples are also given.
It should be noted that this material can be tricky, despite
its apparent simplicity, and one should take care to make precise
statements and proofs as well as to check examples.

If $f(x,y)$ is a polynomial, $p \in \lineinfinity$ is a point
through which the level curves of $f$ pass, and $c \in
\complex$, we say that the pair $(p,c)$ is a {\em regular point at infinity}
for $f(x,y)$ if it satisfies any one of the following equivalent
conditions.  (Otherwise it is a {\em critical point at infinity}.)

\begin{itemize}
\item Condition M (\ref{condition-m}):  There is no jump in the Milnor
number (2.1):  $\nu_{p,c} = 0$.
\item Condition E (\ref{condition-e}): The family of
germs $f(x,y) = tz^d$ at $p$ is equisingular at $t = c$. 
\item Condition F (\ref{condition-f}): The map $f$ is a smooth fiber bundle
near $p$ and the value $c$. 
\item Condition R (\ref{condition-r}): There is a 
resolution $\tilde{f}: M \to \projective^1$ with $\pi: M \to
\projective^2$ and a neighborhood $U$ of $p \in \projective^2$
such that $\{ \tilde{f} = c \} \cap \pi^{-1}(U)$ is smooth and
intersects the exceptional set $ \pi^{-1}(p)$ transversally. 
\item Condition G (\ref{condition-g}): There does not exist 
a sequence of points $ \{ p_k \} \in \complex^2 $
with $ p_k \to p$, $grad \, f(p_k) \to 0$ and $f(p_k) \to c$ as $k \to
\infinity $. 
\end{itemize}

%The equivalence of Conditions M and E is well known, 
%as is the equivalence of Conditions E and F
%(Proposition \ref{E-equivalent-F}).  The equivalence of Conditions E
%and R follows from Zariski's work on equisingularity (Proposition
%\ref{E-equivalent-R}).  Condition G implies Condition F by an
%argument originally due to Broughton (Proposition \ref{G-implies-F}).
%It is shown in \cite{Siersma-Tibar-95} and \cite{Parusinski-P95} that
%Condition M implies Condition G; their argument uses Whitney conditions at
%infinity.

Most of these equivalences are well known; we give either
proofs or references
for proofs in the pages that follow.

There are several new results in this paper.
First, 
we define (\ref{nu-infinity}) an invariant $\nu_{p,
\infinity}$ which measures the number of vanishing cycles at a point
$p$ on the line at infinity for the critical
value infinity, and show that this invariant has
many of the same properties that $\nu_{p,c}$ does for $c \in \complex$.
Secondly, we define $g_{p,c}$ to be the
number of isotopy classes of paths $\alpha : \real \to \complex^2$
such that $\alpha (t) \to p$, $grad \, f(\alpha (t) ) \to 0$ and
$f(\alpha(t)) \to c$ as $t \to + \infinity$.
We use this to give a new proof that Condition M implies Condition G.
In fact, we will show (Proposition \ref{nu-geq-g}) that
$\nu_{p,c} \geq g_{p,c}$.

The work described in this paper started in 1989 when the author supervised a
group of undergraduates in the Mount Holyoke Summer Research Institute
in Mathematics who were working on corresponding problems for real
polynomials.  These results are described in \cite{REU}, with further
results in \cite{Durfee-P95}. 
The work for this paper was carried out
at the Tata Institute, Bombay, Martin-Luther University,
Halle (with support from IREX, the International Research and
Exchanges Board), the University of Nijmegen, Warwick University, the
Massachusetts Institute of Technology and the University of Bordeaux.
The author would like to thank all of them for their hospitality.

Earlier versions of this paper included results on deformations of
critical points at infinity; these will appear elsewhere.

%%%%%%%%%%%%%%%%%%%%%%%%%%%%%

\section{Numerical Invariants}
%\resetassertioncounter

We will use coordinates $(x,y)$ for the complex plane $\complex^2$, and coordinates
$[x,y,z]$ for the projective plane $\complexprojective ^2$.
We let 
$$\clinf = \{ [x,y,z] \in \complexprojective^2: z=0 \}$$ 
be the line at infinity.

We let $d$ be the degree of the polynomial $f(x,y)$.
We let $f_d$ denote the homogeneous term of degree $d$ in $f$.
If $p = [a,b,0] \in \lineinfinity$, we let $d_p$ be the
multiplicity of the factor $(bx-ay)$ in $f_d$.

%We let $m_p$ be the intersection multiplicity of the projective
%completions of the curves 
%$\{ f_x = 0 \}$ and $ \{ f_y = 0 \}$ at $p$.

Suppose that the level sets of $f$ intersect $\lineinfinity$ at $p$.
Let 
$$F_t(x,y,z) = z^df(x/z, y/z) - tz^d$$ 
be the homogenization of the polynomial $f(x,y)-t$,
where $t \in \complex$,
and let $g_{p,t}$ be the local equation of $F_t$ at $p$.  
If $p = [1,0,0]$, then $g_{p.t}$ is given in local coordinates 
$$(u,v) = (y/x, 1/x)$$ 
by
$$g_{p,t}(u,v) = F_t(1,u,v)= v^df(1/v,u/v) - tv^d$$

Note that the multiplicity of $g_{p,t}$ at $(0,0)$ is at most $d_p$.

\begin{definition} 
\label{Milnor-number}
The {\em Milnor number} $\mu_{p,t}$ of $f(x,y)$ at 
$(p,t) \in \lineinfinity \times \complex$ is the Milnor number of the
germ $g_{p,t}$ at $(0,0)$ in the usual sense. 
The {\em generic Milnor number} $\mu_{p,gen}$ is the Milnor number $\mu_{p,t}$ for generic $t$.
The number of {\em vanishing cycles} at $(p,t)$ is 
$$\nu_{p,t} = \mu_{p,t} - \mu_{p,gen}$$
\end{definition}

\begin{example} Let $f(x,y) = y(xy-1)$ and $p = [1,0,0]$.  
Then  $g_{p,t}(u,v) = u^2-uv^2-tv^3$.  
We have $\mu_{p,gen}=2$, $\nu_{p,0}=1$, and all other $\nu_{p,t} = 0$.  
In fact,
for $t \neq 0$, the singularity is of type $A_2$, and for $t=0$, the
singularity is of type $A_3$.  This well-known example is the simplest ``critical point
at infinity".
More generally, if $f(x,y) = y(x^ay-1)$ then for $t \neq 0$, $\mu_{p,t} = a+1$ and there is a singularity of type $A_{a+1}$.  For $t=0$, $\mu_{p,0} = 2a+1$ and there is a singularity of type $A_{2a+1}$.

%\begin{figure}
%\postscript{cfigures/std-crpt-level.eps}{0.7}
%\caption{Level sets of the function $y(xy-1)$}
%\label{std-crpt-level}
%\end{figure}
\end{example}

\begin{example} Let $f(x,y)=x(y^2-1)$ and $p = [1,0,0]$.  
Then $g_{p,t}(u,v) = u^2-v^2-tv^3$.
For all $t$, $\mu_{p,t}=1$; the family is equisingular, and there is
no ``critical point at infinity''.  This is another basic example.

%\begin{figure}
%\postscript{cfigures/std-nocrpt-level.eps}{0.7}
%\caption{Level sets of $x(y^2-1)$}
%\label{std-nocrpt-level}
%\end{figure}
\end{example}

\begin{example} 
\label{two-max-ex}  Here is a more complicated example
(see \cite{REU,Durfee-P95}): 
Let $f(x,y) = (xy^2-y-1)^2 + (y^2-1)^2$.
At $p = [1,0,0]$ we have $\mu_{gen} = 15$, $\nu_{p,1}=2$, $\nu_{p,2} =
1$ and $\nu_{p,c}=0$ for all other $c$.
\end{example}

Next we relate $\nu_{p,t}$ to homological vanishing
cycles.
Fix $p \in \clinf$ and $c \in \complex \cup \{\infinity \}$.
Let  $U \subset \complex^2$ be an open set such that  
the closure in projective space of the set 
$$ \{ (x,y) \in \complex^2 : (x,y) \in U
\mbox{\ and \ }
f(x,y)= t \} $$
is $p$ for $t$ near $c$.
Choose $C > 0$ large.  
We define the {\em Milnor fiber} of $f$ at $(p,c)$ to be
$$\tilde{F}_{p,c} = \overline{ \{ (x,y) \in \complex^2 : (x,y) \in U
\mbox{\ and \ }
|(x,y)| \geq  C
\mbox{\ and \ }
f(x,y)= t \} } $$
where the overbar indicates closure in projective space, and, if $c
\in \complex$, then $t$ is
near, but not equal to, $c$, and if $c = \infinity$, then $t$ is large.

\begin{proposition}
\label{betti}
For $p \in \clinf$ and $c \in \complex$, 
$$\nu_{p,c} = rank \, H_1 (\tilde{F}_{p,c})$$
\end{proposition}

\begin{xproof}
Without loss of generality, we may assume that $p = [1,0,0]$.
The number $\nu_{p,c}$ is the difference of the
Milnor number $\mu_{p,c}$ and the generic Milnor number $\mu_{p,gen}$.
The number $\mu_{p,gen}$ 
is the
Milnor number of $g_{p,t}$ for $t$ near, but not equal to, $c$.
By the usual argument, this difference is 
$rank \, H_1 (\{g_{p,c}= 0 \} \cap B_0)$
where $B_0$ is the small ball for the Milnor number of $g_{p,t}$.
We may replace $\{g_{p,c}= 0 \} \cap B_0$ by
$$F'_{p,c} = \{ (u,v) \in \complex^2 : 
|v| \leq \epsilon' 
\mbox{\ and \ } g_{p,t}(u,v) = 0 \} $$

We may replace $\tilde{F}_{p,c}$ by
$$\tilde{F}'_{p,c} = \overline{ \{ (x,y) \in \complex^2 : (x,y) \in U
\mbox{\ and \ }
|x| \geq  C
\mbox{\ and \ }
f(x,y)= t \} } $$

The change of coordinates 
$x = 1/v$ and $y = u/v$
takes $\tilde{F}_{p,c}$ to $F'_{p,c}$.
\end{xproof}

To define ``vanishing cycles'' for the critical value $c = \infinity$,
we take the above proposition
to be a definition:

\begin{definition}
\label{nu-infinity}
For $p \in \lineinfinity$ we let
$$\nu_{p,\infinity} = rank \, H_1 (\tilde{F}_{p,\infinity})$$
\end{definition}

\begin{remark}
Here is a topological interpretation of the number of vanishing
cycles at infinity: 
Suppose the level curves of the polynomial $f$ of degree $d$ 
intersect $\clinf$ at $k$ points (counted
without multiplicities).
Then $\nu_{p,\infinity} = 0$ for all $p \in \clinf$ if and
only if $\overline{ \{ f(x,y) = t \} }$ for $t$ large is homeomorphic
to a $d$-fold
cover of $\clinf$ branched at $k$ points.
For example, $y(xy-1) = t$ (where $\nu_{p,\infinity} = 0$ for all $p$) is a three-fold cover of $\projective^1$
branched at two points,
but $y^2-x = t$ (where $\nu_{[1,0,0],\infinity} = 1$) is not a two-fold cover of $\projective^1$ branched
at one point.
\end{remark}

Next we describe three ways of computing the number of vanishing
cycles $\nu_{p,c}$.
The first is to compute (perhaps with a computer algebra progam) $\mu_{p,c}$ and $\mu_{p,gen}$ and
subtract. 
The second is by counting nondegerate critical points, as described in
the proposition below. 
This is similar to
computing the usual Milnor number by counting the number
of nondegerate critical points in a Morsification
(see \cite[vol II, p. 31]{AGV}), and the proof is similar.

\begin{proposition}
\label{nu-equals-critical-values}
Let $p \in \lineinfinity$ and $c \in \complex \cup \{ \infinity \}$. 
The number $\nu_{p,c}$ is equal to the
number of critical points (assumed nondegenerate) 
$q \neq (0,0)$ of the function
$g_{p,t}$ such that $q \to (0,0)$ as $t \to c$.
\end{proposition}

%\begin{xproof}
%The Milnor number of $g_{p,c}$ at the point $(0,0)$ is $\mu_{p,gen} +
%\nu_{p,c}$.  Since $g_{p,t}$ is a deformation of $g_{p,c}$, this
%Milnor number equals the sum of the Milnor numbers of $g_{p,t}$ at the
%point $(0,0)$, which is $\mu_{p,gen}$ plus the sum of the Milnor
%numbers at points near $(0,0)$.
%\end{xproof}

\begin{example}
Let $f(x,y) = y(xy-1)$ and $p = [1,0,0]$.  Then $g_{p,t}(u,v) =
u^2-uv^2-tv^3$.  For $t \neq 0$ the function $g_{p,t}$ has a
(degenerate) critical point at $(0,0)$ with critical value 0, and a
nondegenerate critical point at $((9/2) t^2, -3t)$ with critical value
$(27/4)t^4$.  As $t \to 0$ the second critical point approaches
$(0,0)$.  Thus $\nu_{p,0} = 1$.
\end{example}

\begin{example}
Let $f(x,y) = x-y^2$ and $p = [1,0,0]$.  
Then $g_{p,t}(u,v) = v - u^2 - tu^2$.  
The function
$g_{p,t}$ has a single nondegenerate critical point at
$(0,1/(2t))$ with critical value $1/(4t)$.  As $t \to \infinity$ this
critical point approaches $(0,0)$, so $\nu_{p,\infinity} = 1$. 
\end{example}

%\begin{example}
%Let $f(x,y) = x(y^2-1)$ and $p = [1,0,0]$.  Then $g_{p,t}(u,v) =
%u^2-v^2-tv^3$.  The function $g_{p,t}$ has critical points

%\medskip
%\begin{tabular}{|l|l|l|}  \hline
%critical point & critical value & Milnor number \\ \hline 
%$(0,0)$ &$0$ & $1$ \\ \hline 
%$(0, -2/(3t) )$ & $-4/(27t^2)$ & $1$ \\ \hline
%\end{tabular}
%\medskip

%\noindent As $t \to \infinity$ the critical point $(0, -2/(3t) ) \to
%(0,0)$.  Thus $\nu_{p,\infinity} = 1$.
%\end{example}

The next proposition describes the result of computing
$\nu_{p,\infinity}$ by similar methods.

\begin{proposition}
For $p \in \lineinfinity$, 
$$\nu_{p,\infinity} = (d_p-1)(d-1) - \mu_{p,gen}$$
\end{proposition}

\begin{xproof}
Without loss of generality $p = [1,0,0]$.
The intersection multiplicity of the
curves $(g_{p,t})_u$ and $(g_{p,t})_v$ at $(0,0)$ for $t = \infinity$,
where $(u,v)$ are local coordinates at $(0,0)$,
can be computed using the algorithm in \cite{Fulton}, and is
found to be $(d_p-1)(d-1)$.
(To compute the intersection multiplicity at $t = \infinity$, we let
$s = 1/t$ and compute it at $s = 0$.)
For large $t \neq \infinity$, the intersections split into those at
$(0,0)$, the number of which is $\mu_{p,gen}$, and those not at
$(0,0)$, the number of which is
$\nu_{p, \infinity}$.
\end{xproof}

\begin{example}
The polynomial $f(x,y) = y^a + x^{a-2}y+x$ has $\nu_{p,\infinity} = a^2-2a$ at
the point $p = [1,0,0]$, and all other $\nu_{p,c} = 0$.
\end{example}

Finally, $\nu_{p,t}$ can computed a third way by using polar curves,
as described below. (See \cite[1.6, 1.8]{Ha-Nguyen-89}.)
This method also shows that some vanishing cycles are easy to ``see'' from
a contour plot, since they are where the level
curves of the polynomial have a vertical tangent.

\begin{proposition}
\label{polar-curves}
Suppose $p = [1,0,0]$, $c \in \complex \cup \{\infinity\}$ and the level sets of $f$ pass through $p$.  
Then $\nu_{p,c}$ 
is the number of points of intersection 
$q \in \complex^2$ (assumed transverse) of the curves $f =
t$ and $f_y = 0$ in $\complex^2$ such that $q \to p$ as $t \to c$.
\end{proposition}

\begin{xproof}
The set $F'_{p,c}$ from the proof of Proposition \ref{betti}
%$\overline{ \{ f(x,y) = t \} } \cap B_p$ 
is a connected branched cover of the disk $|v| \leq \epsilon'$ in the
$uv$-plane. 
Two sheets come together at each branch point, and all the sheets come
together over $p$.
The result follows from Hurwitz's
formula.
\end{xproof}

\begin{example}
If $f(x,y) = x(y^2-1)$, the curves $f=t$ and $f_y =0$
intersect at $(-t,0)$.  As $t \to \infinity$, the intersection point
$(-t,0) \to [1,0,0]$ and $f(t,0) \to \infinity$.  Thus
$\nu_{[1,0,0],\infinity} = 1$.  All other $\nu_{[1,0,0],c} = 0$.
\end{example}

%%%%%%%%%%%%%%%%%%%%%%%%%%%%%%%%%%%%%%%%%%%%%%%%%%%%

Next we give three definitions of ``critical
point at infinity''.

\begin{definition}
\label{condition-m} 
The polynomial $f(x,y)$ satisfies Condition M at the point 
$p \in \lineinfinity$ and $c \in \cinf$
 if $\nu_{p,c} = 0$.
\end{definition}

\begin{definition} 
\label{condition-e}
The polynomial $f(x,y)$ satisfies Condition E at the point 
$(p,c) \in \lineinfinity \times  \complex$ if the family of
germs $g_{p,t}$ at (0,0) is equisingular at $t = c$.
\end{definition}

A proof that Condition M for $c \in \complex$ is equivalent to Condition E 
may be found at the end of
\cite{Le-Ramanujam}).

There are various equivalent ways of specifying equisingularity; see
for instance the papers by Zariski in volume IV of \cite{Zariski-works}.  
One that will be useful for us is the following:
The family of germs $g_{p,t}$ is equisingular if the germs  
$g_{p,t}= 0$ at (0,0) form a fiber bundle near $t = c$.

%More precisely, suppose $p = [1,0,0]$.
%They form a fiber bundle if
%there is an $\epsilon > 0$ and a $\delta > 0$ such that, letting
%$$D = \{ t \in \complex : |t - c | < \delta \} $$
%and
%$$M = \{ (u,v,t) \in \complex^2 \times \complex : 
%|(u,v)| \leq \epsilon 
%\mbox{\ and \ } t \in D 
%\mbox{\ and \ } g_{p,t}(u,v) = 0 \} $$
%then the projection to the third coordinate
%$$\pi: M \to D $$
%is a fiber bundle.

\begin{definition}
\label{condition-f}
The polynomial $f(x,y)$ satisfies Condition F at a point $(p,c) \in
\lineinfinity \times \complex$ if the map $f$ is a smooth fiber bundle
near $p$ and the value $c$.
(More precisely, a polynomial satisfies Condition F if 
there is a $U \subset \complex^2$ with $p$ in the closure of $U$ in
projective space and $C > 0$ and $\beta > 0$ such
that, letting
$$B = \{ t \in \complex : |t - c | \leq \beta \}$$
and 
$$N = \{ (x,y) \in \complex^2 : (x,y) \in U
\mbox{\ and \ }
|(x,y)| \geq  C
\mbox{\ and \ }
f(x,y) \in B \} $$
then
$$f: N \to B $$
is a smooth fiber bundle.)
\end{definition}

\begin{proposition}
\label{E-equivalent-F}
A polynomial $f(x,y)$ satisfies Condition E at a point $(p,c) \in
\lineinfinity \times \complex$ 
if and only if it satisfies Condition F at that point.
\end{proposition}

\begin{xproof}
The proof is straight-forward, and just involves replacing the
``spherical'' Milnor fiber by one in a ``box'':
Without loss of generality, we may assume that $p = [1,0,0]$.
We may replace Condition E by the following:
There is an $\epsilon' > 0$ and a $\delta' > 0$ such that, letting
$$D' = \{ t \in \complex : |t - c | < \delta' \} $$
and
$$M' = \{ (u,v,t) \in \complex^2 \times \complex : 
|v| \leq \epsilon' 
\mbox{\ and \ } t \in D' 
\mbox{\ and \ } g_{p,t}(u,v) = 0 \} $$
then the restriction of the projection to the third coordinate
$$\pi: M' \to D' $$
is a fiber bundle.
We may do this since the germs $g_{p,t}(u,v) = 0$ never have $v = 0$
as a component.

We may also replace Condition F by the following:
There is a $U' \subset \complex^2$  
with $p$ in the closure of $U'$ in projective space
and $C' > 0$ and $\beta' > 0$ such
that, letting
$$B' = \{ t \in \complex : |t - c | \leq \beta' \}$$
and 
$$N' = \{ (x,y) \in \complex^2 : (x,y) \in U'
\mbox{\ and \ }
|x| \geq  C'
\mbox{\ and \ }
f(x,y) \in B' \} $$
then
$$f: N' \to B' $$
is a smooth fiber bundle.
 
The change of coordinates 
$x = 1/v$ and $y = u/v$
takes $f(x,y) = t$ to $g_{p,t}(u,v)= 0$ and
$N'$ to $M'$.
\end{xproof}

%%%%%%%%%%%%%%%%%%%%%%%%%%%%%%%%

\section{Resolutions}
%\resetassertioncounter
\label{sec-resolutions}

The polynomial $$f: \complex^2 \to \complex $$ extends to a map
$$\hat{f} : \projective^2 \to \projective $$ which is undefined at a
finite number of points on the line at infinity $\lineinfinity$.
By blowing up these points one gets a manifold $M$ and a map $$ \pi : M \to \projective^2 $$ such
that the map $$ \tilde{f} : M \to \projective $$ which is the lift of
$\hat{f}$ is everywhere defined.  We call the map $ \tilde{f}$ a
{\em resolution of $f$}.  
Some interesting results on the structure of resolutions are announced
in \cite[Theorems 2, 3, 4]{Le-Weber-95}.

\begin{figure}
\postscript{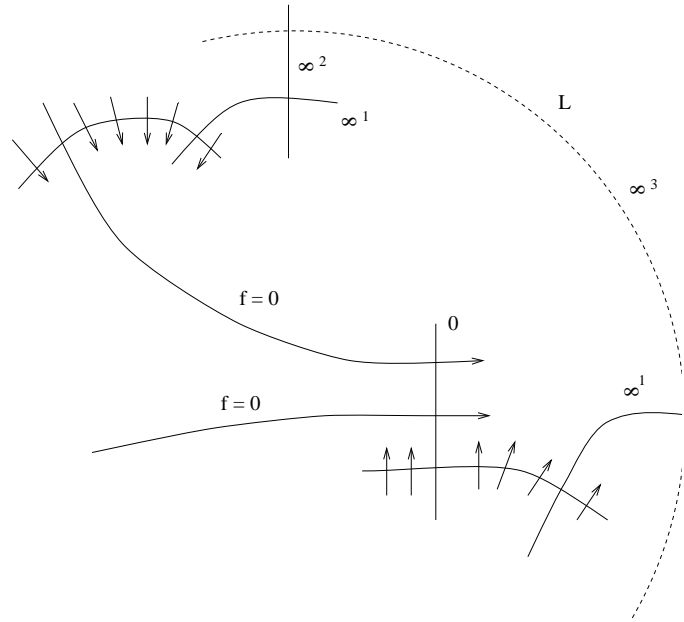}{0.7}
\caption{Resolution of $y(xy-1)$}
\label{std-crpt-res-g}
\end{figure}

For example, a resolution (the minimal resolution) of 
$y(xy-1)$ is given in Figure \ref{std-crpt-res-g}.
The symbol $c^{m}$ next to a divisor means that at each smooth
point of the divisor there are local coordinates $(z,w)$ in a
neighborhood of the point such that the divisor is $z=0$ and
$\tilde{f}(z,w) =(z-c)^m$. 
The proper transform of level curves of $f$ have arrowheads on them;
the exceptional sets do not.

\begin{figure}
\postscript{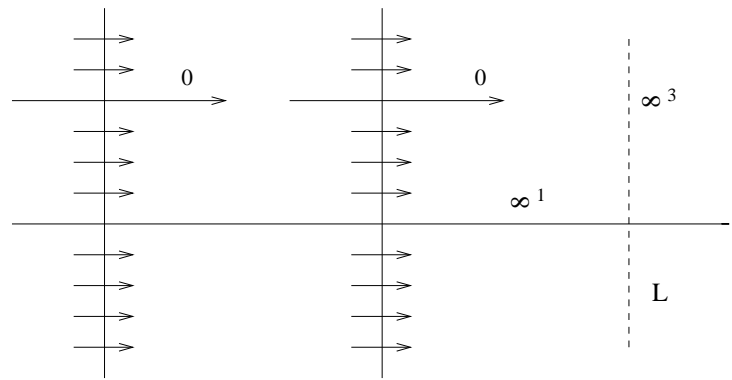}{0.7}
\caption{Resolution of $x(y^2-1)$ at $[1,0,0]$}
\label{std-nocrpt-res}
\end{figure}

\begin{figure}
\postscript{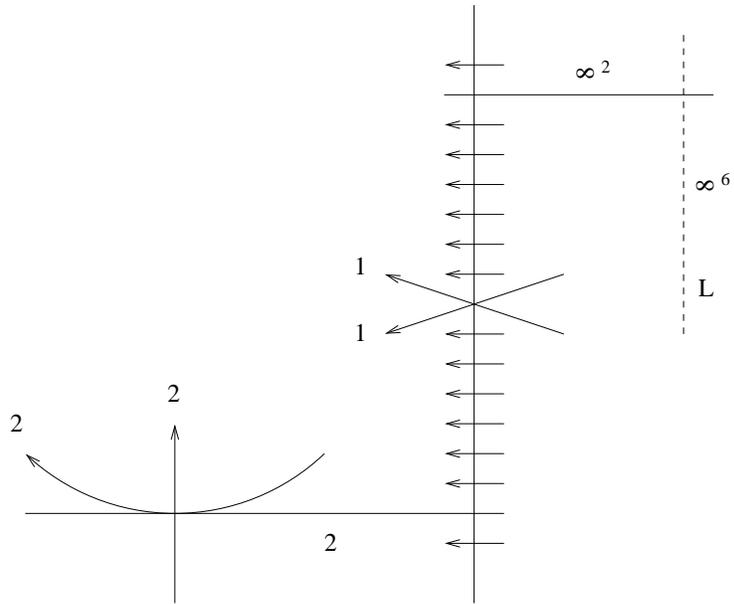}{0.7}
\caption{Resolution of $(xy^2-y-1)^2 + (y^2-1)^2$ at $[1,0,0]$}
\label{two-max-res}
\end{figure}

Resolution are easy compute.
For example, starting with $f(x,y) = y(xy-1)$ which we wish to resolve
near $[1,0,0]$, the function in local coordinates at $[1,0,0]$ is
$u(u-v^2)/v^3$, and we blow up in the standard fashion until it is
everywhere defined.
More examples are shown in Figures \ref{std-nocrpt-res} and
\ref{two-max-res}.

%\begin{figure}
%\postscript{cfigures/max-min-res.eps}{0.7}
%\caption{Resolution of $y^5 + x^2y^3 - y$ at $[1,0,0]$}
%\label{max-min-res}
%\end{figure}

Next we give a condition for ``regular point at infinity'' in terms of
a resolution. (See also \cite[Theorem 5]{Le-Weber-95}.)

\begin{definition}
\label{condition-r}
The polynomial $f(x,y)$ satisfies Condition R at a point $(p,c) \in
\lineinfinity \times \complex$ if there is a 
resolution $\tilde{f}: M \to \projective^1$ with $\pi: M \to
\projective^2$ and a neighborhood $U$ of $p \in \projective^2$
such that $\{ \tilde{f} = c \} \cap \pi^{-1}(U)$ is smooth and
intersects the exceptional set $ \pi^{-1}(p)$ transversally.
\end{definition}

\begin{example}
\label{Krasinski}
Let $f(x,y) =y-(xy-1)^2$ near $[1,0,0]$ (See \cite{Krasinski}).  In
this example the level curve of the function $\tilde{f} = 0$ is
smooth, but it does not intersect the exceptional divisor
transversally; see Figure \ref{kras-res}.
Hence $(p,c) = ([1,0,0],0)$ does not satisfy Condition R.
(Here $\nu_{[1,0,0],0} = 1$.)
\begin{figure}
\postscript{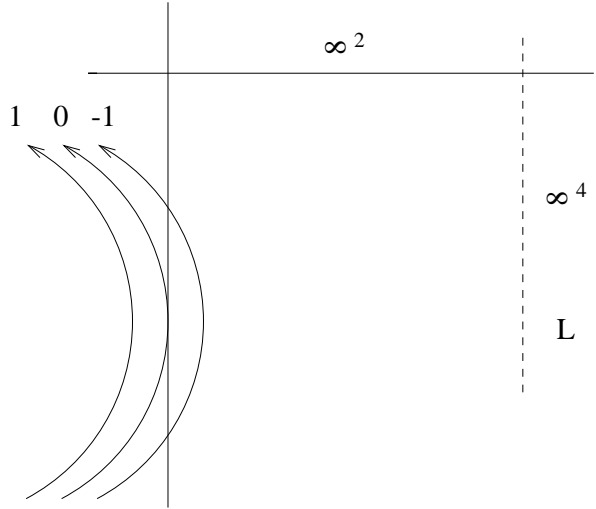}{0.7}
\caption{Resolution of $y-(xy-1)^2$ at $[1,0,0]$}
\label{kras-res}
\end{figure}
\end{example}

\begin{proposition}
\label{E-equivalent-R}
A point $p \in \lineinfinity$ and a value $c \in \complex$ for a
polynomial $f(x,y)$ satisfies Condition E if and only if it satisfies
Condition R.
\end{proposition}

\begin{xproof}
Suppose $(p,c)$ satisfies Condition E.
Let $U$ be a neighborhood of $p$ in $\complexprojective ^2$ containing
no critical points of $f$ in $\complex ^2$ or points on $\lineinfinity$ though
which the level curves of $f$ pass.
Find a resolution $\tilde{f}$ of $f$.
By further blowing up (if necessary), we may assume that 
$\tilde{f}^{-1}(c)$
is a divisor with normal crossings transversally intersecting the
exceptional set where $\tilde{f}$ is not constant.
Equisingularity in
the form of Zariski's (b)-equivalence \cite[p. 513]{Zariski-65}
implies that the functions $g_{p,t}$ for $t$ near $c$ have the same
resolution as the function $g_{p,c}$. 
This can only happen if $\{\tilde{f}^{-1}(t) \cap \pi^{-1}(U) \}$ is
smooth and transversally intersects the exceptional set $\pi ^{-1}(p)$.
Thus $p$ and $c$ satisfy Condition R.

Conversely, if $p$ and $c$ satisfy Condition R, then the resolutions
of  $ \{ g_{p,t} = 0 \}$ for $t$ near $c$ 
are (b)-equivalent and hence equisingular.
\end{xproof}

%%%%%%%%%%%%%%%%%%%%%%%%%%%%

%%%%%%%%%%%%%%%%%%%%%%%%%%%%%

\section{The Gradient}
\resetassertioncounter

If $f$ is a complex polynomial, we define $grad \, f$ as in
\cite{Milnor} to be the complex conjugate of the vector of partial
derivatives.

Of course $p \in \complex^2$ is a regular point for a function $f$
with regular value $c \in \complex$ if $f(p) = c$ and $grad \, f(p) \neq 0$.  An
equivalent definition would be to say that there is no sequence of
points $ \{ p_k \} $ with $ p_k \to p$, $grad \, f(p_k) \to 0$ and $f(p_k) \to c$ as $k \to
\infinity $.  We can now imitate this definition for $p \in \lineinfinity$
as follows:

\begin{definition}
\label{condition-g}
The polynomial $f(x,y)$ satisfies Condition G at a point $p \in
\lineinfinity$ and $c \in \cinf$ if
there does not exist a sequence of points $ \{ p_k \} \in \complex^2 $
with $ p_k \to p$, $grad \, f(p_k) \to 0$ and $f(p_k) \to c$ as $k \to
\infinity $.
\end{definition}

If $(p,c)$ does not satisfy Condition G, then a
version of Milnor's curve selection lemma (see for instance
\cite[Lemma 3.1]{Ha-P91} or \cite[Lemma 2]{Nemethi-Zaharia-92}) 
implies that the sequence of points can be
replaced by a curve:

\begin{lemma}
\label{curveselectionlemma}
If $(p,c)$ does not satisfy Condition G, then there is a smooth real algebraic curve $\alpha :
\real^+ \to \complex^2$ such that $\alpha(t) \to p$, \ $grad \,
f(\alpha(t)) \to 0$ and $f(\alpha(t)) \to c$ as $t \to + \infinity$.
\end{lemma}

By ``real algebraic curve'' we mean that the image of $\alpha$ in
$\complex^2$ is contained in an irreducible component of the zero locus
of a real polynomial.

\begin{example}
Let $f(x,y) = y(xy-1)$.
%There is a critical point at $[1,0,0]$ with critical value $0$.
Let $\alpha(t) = (t, 1/(2t))$.  As $t \to + \infinity$, $\alpha(t) \to
[1,0,0]$, the gradient of $f$ goes to $0$ and the value of the
function approaches $0$.
\end{example}

\begin{example}
(c.f. Example \ref{two-max-ex}.)
Let $f(x,y) = (xy^2-y-1)^2 + (y^2-1)^2$.
%There is a critical point at $[1,0,0]$ with two critical values, $-1$ and $-2$.
Let $\alpha(t) = (t+t^2, \pm 1/t)$.  As $t \to + \infinity$, $\alpha(t)
\to [1,0,0]$, the gradient of $f$ goes to $0$ and the function
approaches the value $1$.  If $\beta(t) = (t/2,
1/t)$, then as $t \to + \infinity$, $\beta(t) \to [1,0,0]$, the
gradient of $f$ goes to $0$ and the function approaches the value
$2$.  (These paths were found by Ian Robertson in the Mount Holyoke
REU program in the summer of 1992.)
\end{example}

%\begin{example}
%Let $f(x,y) = y - (xy-1)^2$.
%Let $\alpha(t) = (t, (2t+1)/(2t^2))$.
%As $t \to \infinity$, $\alpha(t) \to [1,0,0]$, $grad \, f(\alpha(t))
%\to 0$ and $f(\alpha(t)) \to 0$.
%\end{example}

\begin{example}
If $f(x,y) = x^2y + xy^2 + x^5y^3 + x^3y^5$ 
and $q \to [1,0,0]$ along the curve $y^2x^3 = -1/3$,
then $grad \, f(q) \to 0$ and $f(q) \to \infinity$.
Here $v_{p,\infinity} = 1$. 
This polynomial is ``quasi-tame'' but not ``tame'' 
\cite{Nemethi-Zaharia-92}.
It would be interesting to find more examples like this.
\end{example}

The following proposition is well-known.
It was first proved in the global case by Broughton
\cite{Broughton-88}; see also \cite{Nemethi-Zaharia-90}, proof of
Theorem 1, and \cite{Siersma-Tibar-95}, proof of Proposition 5.5.
The idea of the proof is to use integral curves of the vector field 
$grad \, f / |grad \, f |^2$ to identify the fibers.

\begin{proposition}
\label{G-implies-F}
If a polynomial $f(x,y)$ satisfies Condition G at $(p,c) \in \lineinfinity \times \complex$, then it satisfies Condition
F at this point.
\end{proposition}

Next we will show that Condition M implies Condition G;
this has been shown in \cite{Ha-90, Ha-P91, Siersma-Tibar-95,
Parusinski-P95}.
Here we will prove a stronger result by
different methods.

\begin{definition}
\label{def-gpc1}
For $p \in \lineinfinity$ and $c \in \cinf$, let $g_{p,c}$ be the
number of isotopy classes of smooth real algebraic curves $\alpha : \real \to \complex^2$
such that $\alpha (t) \to p$, $grad \, f(\alpha (t) ) \to 0$ and
$f(\alpha(t)) \to c$ as $t \to + \infinity$.
\end{definition}

\begin{example} 
If $f(x,y) = y^5+x^2y^3-y$ and $p = [1,0,0]$, then
$\nu_{p,0} = 2$.  There are two isotopy classes of curves approaching
$p$ along which $grad \, f$ goes to zero, namely the ones containing
the two branches of the curve $f_y = 0$ at $p$.  Hence $g_{p,0} =2$.
(This example is from \cite{REU}.)
\end{example}

Clearly $(p,c) \in \lineinfinity \times (\cinf)$ satisfies Condition G
if and only if $g_{p,c} = 0$.  
Now let $\pi : M \to \projective^2$ be a resolution of $f$, $f_x$, and
$f_y$ (so that $\tilde{f}$, $\widetilde{(f_x)}$ and
$\widetilde{(f_y)}$ are defined on $M$), and let
$$G_{p,c} = \{ q \in M : \pi(q) = p,\ \tilde{f}(q) = c,\ \widetilde{(f_x)}(q)
= 0 \mbox{ \ and \ } \widetilde{(f_y)}(q) = 0 \} $$

\begin{definition}
\label{def-gpc2}  
For $p \in \lineinfinity$ and $c \in \cinf$, let 
$\tilde{g}_{p,c}$ be the number of connected components of $G_{p,c}$.
\end{definition}

The number $\tilde{g}_{p,c}$ is independent of the resolution by the
usual argument.  

\begin{example}
In the minimal resolution of $f(x,y) = y(xy-1)$ at $p = [1,0,0]$ (Figure
\ref{std-crpt-res-g}), the functions $f_x$ and $f_y$ are defined.
The zero locus of the lift of $f_x$ contains the exceptional set where
the lift of $f$ is zero, and
the zero locus of the lift of $f_y$ intersects this set transversally.
Thus $G_{p,0}$ consists of a single point, and $\tilde{g}_{p,0} = 1$.
If $f(x,y) = y^5+x^2y^3-y$ and $p = [1,0,0]$, one finds similarly that
$G_{p,0}$ consists of two points.
\end{example}

The two definitions are equivalent by the following
proposition.

\begin{proposition}  For $p \in \lineinfinity$ and $c \in \cinf$,
$g_{p,c} = \tilde{g}_{p,c}$.
\end{proposition}

\begin{xproof}
Let $\pi: M \to \projective^2$ be a resolution of $f$, $f_x$ and $f_y$.
We will show that there is a one-one correspondence between isotopy
classes of curves satisfying (\ref{def-gpc1}) and connected
components of $G_{p,c}$.
Suppose that $\alpha : \real^+ \to \complex^2$ is a smooth real
algebraic curve satisfying the conditions of (\ref{def-gpc1}).
Since $\alpha$ is real algebraic, it lifts to to a map 
$\tilde{\alpha}: \real^+ \cup \{ \infinity \} \to M$
with $\tilde{\alpha}(\infinity) \in \pi^{-1}(\lineinfinity)$.
Let $q = \tilde{\alpha}(\infinity)$.
Then $\tilde{f}(q) = c$ and $\widetilde{(f_x)}(q) = 0$ and
$\widetilde{(f_y)}(q) = 0$.
Thus $q \in G_{p,c}$.

If $\alpha_0$ is isotopic to $\alpha_1$ through curves $\alpha_t$
satisfying (\ref{def-gpc1}), then the curves $\alpha_t$ lift to $M$
and are isotopic.  In particular, $\tilde{\alpha_0}(\infinity)$ and
$\tilde{\alpha_1}(\infinity)$ are in the same connected component of
$G_{p,c}$.

For each $q \in G_{p,c}$ there is an algebraic curve $\tilde{\alpha}:
\real^+ \cup \{ \infinity \} \to M$ with $\tilde{\alpha}(\infinity) =
q$ and $\tilde{\alpha}(R^+) \subset \pi^{-1}(\complex^2)$.  Let
$\alpha = \pi \circ \tilde{\alpha}: \real^+ \to \complex^2$.  
Then $\alpha$ satisfies the conditions of (\ref{def-gpc1}).
%$\alpha(t) \to p$, \ $grad \, f(\alpha(t)) \to 0$ and $f(\alpha(t))
%\to c$ as $t \to + \infinity$.  
If $q_0, q_1 \in G_{p,c}$, then there
are two such curves $\tilde{\alpha}_0, \tilde{\alpha}_1$.  If $q_0$
and $q_1$ are in the same connected component of $G_{p,c}$, then
$\tilde{\alpha}_0$ is isotopic to $\tilde{\alpha}_1$ through a family
of such curves $\tilde{\alpha}_t$.  Hence $\alpha_0$ is isotopic to
$\alpha_1$ through curves satisfying (\ref{def-gpc1}).
If $q_0$ and $q_1$ are in different connected components,
then $\alpha_0$ is not isotopic to $\alpha_1$ through curves
satisfying (\ref{def-gpc1}).  
Thus $g_{p,c} = \tilde{g}_{p,c}$.
\end{xproof}

%One must be careful with the blowing-up process:
%For example, the polynomial $y - (xy-1)^2$ at the point $[1,0,0]$ 
%becomes defined after two blow-ups but the partial derivative with
%respect to $y$ is not defined, but needs another blow-up.

\begin{proposition}
\label{nu-geq-g}
For $p \in \lineinfinity$ and $c \in \cinf$,
$$\nu_{p,c} \geq g_{p,c}$$
\end{proposition}

\begin{xproof}
We will show that $\nu_{p,c} \geq \tilde{g}_{p,c}$ and will
use Proposition \ref{polar-curves} to compute $\nu_{p,c}$.
We may assume without loss of generality that $p = [1,0,0]$.
Pick a connected component $G'$ of $G_{p,c}$.
Let $t$ be near $c$.
We will show that $f = t$ intersects $f_y = 0$ in
$\complex^2$ near $G'$.

There is a $q \in G'$ and a component $C$ of $f_y = 0 $ in $\complex
^2$ such that
$q$ is in the closure of $C$ in $M$:  
We have that $\tilde{f_y} = 0$ on $G'$.
Blow down $G'$ to a point $q'$, and let $E'$ be the image of
$\pi^{-1}(p)$.
Then the lift of $f_y$ is not constant on $E'$ near $q'$, so there is
a component of $f_y= 0$ passing through $q'$.  
Lift this component back to $M$.

Next, $f$ is not constant on $C$:  If it were, then the gradient
vector of $f$ would be horizontal, so $C$ would be of the form $x =
const$, and $p$ would not be in the closure of $C$.

Thus $f = t$ intersects $C$ near $q$ for small $\epsilon
\neq 0$, and the intersection points are in $\complex^2$.
\end{xproof}

\begin{remark}
The inequality of the proposition can be strict, as it is for the polynomial
$y(x^2y-1)$ at $p = [1,0,0]$ and $c = 0$, where $\nu_{p,c} = 2$ and $g_{p,c} = 1$.
\end{remark}

\begin{corollary}
If $p \in \lineinfinity$ and $c \in \cinf$ satisfy Condition M, then
they satisfy Condition G.
\end{corollary}

\begin{remark}
The converse to this corollary is not true for $c = \infinity$:
For example, the polynomial $x(y^2-1)$ has a gradient whose
magnitude is bounded below for large $x$
and hence satisfies Condition G at $p = [1,0,0]$, yet
$\nu_{p,\infinity} = 1$.
\end{remark}

%%%%%%%%%%%%%%%%%%%%%%%%%%%%%%%%%%%

%%%%%%%%%%%%%%%%%%%%%%%%%

%%%%%%%%%%%%%%%%%%%%%%%%%%%%%%%%%%%%

\bibliographystyle{alpha}

\begin{thebibliography}{DKM{\etalchar{+}}93}

\bibitem[AGZV85]{AGV}
V.~I. Arnold, S.~M. Gusein-Zade, and A.~N. Varchenko.
\newblock {\em Singularities of Differentiable Maps}.
\newblock Birkhauser, Boston, 1985.

\bibitem[BCND96]{Bartolo-Cassou-Dimca-P96}
E.~Artal Bartolo, Pierrette Cassou-Nogues, and Alexandru Dimca.
\newblock Sur la topologie des polynomes a deux variables complexes.
\newblock Preprint, University of Bordeaux, 1996.

\bibitem[Bro83]{Broughton-83}
S.~A. Broughton.
\newblock On the topology of polynomial hypersurfaces.
\newblock In P.~Orlik, editor, {\em Proceedings of Symposia in Pure
  Mathematics, Vol. 40}, pages 167--178, Providence RI, 1983. American
  Mathematical Society.

\bibitem[Bro88]{Broughton-88}
S.~A. Broughton.
\newblock Milnor numbers and the topology of polynomial hypersurfaces.
\newblock {\em Inventiones Math.}, 92:217--241, 1988.

\bibitem[CN96]{Cassou-P96}
Pierrette Cassou-Nogues.
\newblock Sur la generalisation d'un theoreme de {K}ouchnirenko.
\newblock Preprint, University of Bordeaux, 1996.

\bibitem[CNH94]{Cassou-Ha-P92}
Pierrette Cassou-Nogues and Huy~Vui Ha.
\newblock Theoreme de {K}uiper-{K}uo-{B}ochnak-{L}ojasiewicz a l'infini.
\newblock Preprint, Universite Bordeaux, 1994.

\bibitem[CNH95]{Cassou-Ha-95}
Pierrette Cassou-Nogues and Huy~Vui Ha.
\newblock Sur le nombre de {L}ojasiewicz a l'infiniti d'un polynome.
\newblock {\em Annales Polonici Math.}, 62:23--44, 1995.

\bibitem[DKM{\etalchar{+}}93]{REU}
A.~Durfee, N.~Kronenfeld, H.~Munson, J.~Roy, and I.~Westby.
\newblock Counting critical points of a real polynomial in two variables.
\newblock {\em American Math Monthly}, 100:255--271, 1993.

\bibitem[Dur]{Durfee-P95}
Alan~H. Durfee.
\newblock The index of $grad \, f(x,y)$.
\newblock Duke alg-geom preprint 9506002. To appear (in revised form), {\it
  Topology}.

\bibitem[Ful69]{Fulton}
William Fulton.
\newblock {\em Algebraic Curves}.
\newblock Benjamin, New York, 1969.

\bibitem[Ha89]{Ha-89-2}
Huy~Vui Ha.
\newblock Sur la fibration globale des polynomes de deux variables complexes.
\newblock {\em C. R. Acad. Sci Paris}, 309:231--234, 1989.

\bibitem[Ha90]{Ha-90}
Huy~Vui Ha.
\newblock Nombres de {L}ojasiewicz et singularites a l'infini des polynomes de
  deux variables complexes.
\newblock {\em C. R. Acad. Sci Paris}, 311:429--432, 1990.

\bibitem[Ha91a]{Ha-P91}
Huy~Vui Ha.
\newblock On the irregular at infinity algebraic plane curves.
\newblock Preprint, Mathematical Institute, Hanoi, 1991.

\bibitem[Ha91b]{Ha-91}
Huy~Vui Ha.
\newblock Sur l'irregularite du diagramme splice pour l'entrelacement a
  l'infini des courbes planes.
\newblock {\em C. R. Acad. Sci Paris}, 313:277--280, 1991.

\bibitem[Ha94]{Ha-94}
Huy~Vui Ha.
\newblock A version at infinity of the {K}uiper-{K}uo theorem.
\newblock {\em Acta Math Vietnamica}, 19:3--12, 1994.

\bibitem[HL84]{Ha-Le-84}
Huy~Vui Ha and Dung~Trang Le.
\newblock Sur la topologie des polyn\^omes complexes.
\newblock {\em Acta Math Vietnamica}, 9:21--32, 1984.

\bibitem[HN89]{Ha-Nguyen-89}
Huy~Vui Ha and Le~Anh Nguyen.
\newblock Le comportement geometrique a l'infini des polynomes de deux
  variables complexes.
\newblock {\em C. R. Acad. Sci Paris}, 309:183--186, 1989.

\bibitem[Kra91]{Krasinski}
T.~Krasinski.
\newblock On branches at infinity of a pencil of polynomials in two complex
  variables.
\newblock {\em Ann. Polon. Math.}, 55:213--220, 1991.

\bibitem[LO93]{Le-Oka-P93}
Van~Thanh Le and Mutsuo Oka.
\newblock Estimation of the number of the critical values at infinity of a
  polynomial function $f: C^2 \to c$.
\newblock Preprint, Tokyo Institute of Technology, 1993.

\bibitem[LR76]{Le-Ramanujam}
Dung~Trang Le and C.~P. Ramanujam.
\newblock The invariance of {M}ilnor number implies the invariance of the
  topological type.
\newblock {\em Amer. J. Math}, 98:67--78, 1976.

\bibitem[LW95]{Le-Weber-95}
Dung~Trang Le and Claude Weber.
\newblock Polynomes a fibres rationnelles et conjecture jacobienne a 2
  variables.
\newblock {\em C. R. Acad. Sci Paris}, 320:581--584, 1995.

\bibitem[LW96]{Le-Weber-P96}
Dung~Trang Le and Claude Weber.
\newblock Equisingularite dans les pinceaux de germes de courbes planes et
  $c^0$-suffisance.
\newblock Preprint, University of Geneva, 1996.

\bibitem[Mil68]{Milnor}
John Milnor.
\newblock {\em Singular points of complex hypersurfaces}.
\newblock Princeton University Press, Princeton, 1968.

\bibitem[Neu89]{Neumann-89}
Walter~D. Neumann.
\newblock Complex algebraic plane curves via their links at infinity.
\newblock {\em Invent. math.}, 98:445--489, 1989.

\bibitem[NZ90]{Nemethi-Zaharia-90}
A.~N\'emethi and A.~Zaharia.
\newblock On the bifurcation set of a polynomial function and {N}ewton
  boundary.
\newblock {\em Publ. Math {RIMS} {K}yoto}, 26:681--689, 1990.

\bibitem[NZ92]{Nemethi-Zaharia-92}
A.~N\'emethi and A.~Zaharia.
\newblock Milnor fibration at infinity.
\newblock {\em Indag Math}, 3:323--335, 1992.

\bibitem[Par95]{Parusinski-P95}
Adam Parusinski.
\newblock A note on singularities at infinity of complex polynomials.
\newblock Preprint no. 426, University of {N}ice-{S}ophia-{A}ntipolis, 1995.

\bibitem[Pha83]{Pham-83}
Frederic Pham.
\newblock Vanishing homologies and the $n$ variable saddlepoint method.
\newblock In P.~Orlik, editor, {\em Proceedings of Symposia in Pure
  Mathematics, Vol. 40}, pages 319--333, Providence RI, 1983. American
  Mathematical Society.

\bibitem[ST95]{Siersma-Tibar-95}
Dirk Siersma and Mihai Tibar.
\newblock Singularities at infinity and their vanishing cycles.
\newblock {\em Duke Math. J.}, 80:771--783, 1995.

\bibitem[Suz74]{Suzuki-74}
Masakazu Suzuki.
\newblock Proprietes topologiques des polynomes de deux variables complexes, et
  automorphismes algebriques de l'espace $c^2$.
\newblock {\em J. Math.Soc. Japan}, 26:141--157, 1974.

\bibitem[Tib96]{Tibar-P96}
Mihai Tibar.
\newblock Topology at infinity of polynomial mappings and {T}hom regularity
  condition.
\newblock Preprint, University of Angers, 1996.

\bibitem[Zar65]{Zariski-65}
O.~Zariski.
\newblock Studies in equisingularity i.
\newblock {\em Amer. J. Math.}, 87:507--536, 1965.

\bibitem[Zar79]{Zariski-works}
Oscar Zariski.
\newblock {\em Collected Papers}.
\newblock The MIT Press, Cambridge, 1979.

\end{thebibliography}

\newcommand{\etalchar}[1]{$^{#1}$}

\bigskip

\noindent Department of Mathematics

\noindent Mount Holyoke College

\noindent South Hadley, MA 01075

\medskip

\noindent email:  adurfee@mtholyoke.edu

\end{document}